\documentclass[12pt,reqno]{amsart}
\usepackage{amsmath}
\usepackage{fullpage}

\theoremstyle{plain}
\newtheorem{thm}{Theorem}[section]
\newtheorem{lem}[thm]{Lemma}
\newtheorem{pro}[thm]{Proposition}
\theoremstyle{definition}
\newtheorem{defn}[thm]{Definition}
\newtheorem{ex}[thm]{Example}

\newtheorem{cor}[thm]{Corollary}
\numberwithin{equation}{section}
\newcommand{\R}{\mathbb{R}}

\newcommand{\F}{\mathcal{F}}

\begin{document}
\title[A Birkhoff--Kellogg type result with applications]{An affine Birkhoff--Kellogg type result in cones with applications to functional differential equations}\date{}


\subjclass[2020]{Primary 47H10, secondary 34K10, 34B10, 34B18}%
\keywords{Fixed point index, cone, Birkhoff--Kellogg type result, retarded functional differential equation, functional boundary condition.}%

\author[A. Calamai]{Alessandro Calamai}
\address{Alessandro Calamai, 
Dipartimento di Ingegneria Civile, Edile e Architettura,
Universit\`{a} Politecnica delle Marche
Via Brecce Bianche
I-60131 Ancona, Italy}%
\email{calamai@dipmat.univpm.it}%

\author[G. Infante]{Gennaro Infante}
\address{Gennaro Infante, Dipartimento di Matematica e Informatica, Universit\`{a} della
Calabria, 87036 Arcavacata di Rende, Cosenza, Italy}%
\email{gennaro.infante@unical.it}%

\begin{abstract}
In this short note we prove, by means of classical fixed point index, an affine version of a Birkhoff--Kellogg type theorem in cones. We apply our result to discuss the solvability of a class of boundary value problems for functional differential equations subject to functional boundary conditions. We illustrate our theoretical results in an example.
\end{abstract}

\maketitle

\centerline{\it
Dedicated to Professor Jean Mawhin on the occasion of his eightieth birthday.
}

\section{Introduction}
 The celebrated Birkhoff-Kellogg invariant-direction Theorem~\cite{B-K-1922} is a widely studied and applied tool of nonlinear functional analysis, also in view of its applicability to eigenvalue problems for ODEs and PDEs (see for example the book~\cite{ADV} and the recent papers~\cite{gi-BK, giaml22}). Among the various extensions of the  invariant-direction Theorem, one of them is set in the framework of cones and is due to Krasnosel'ski\u{i} and Lady\v{z}enski\u{\i}~\cite{Kra-Lady}. 
 Before we state this latter result let us 
recall that a cone $K$ of a real Banach space $(X,\| \, \|)$ is a closed set with $K+K\subset K$, $\mu K\subset K$ for all $\mu\ge 0$ and $K\cap(-K)=\{0\}$.  The Birkhoff-Kellogg type theorem of Krasnosel'ski\u{i} and Lady\v{z}enski\u{\i} reads as follows.
\begin{thm}\label{B-K}~\cite[Theorem~2.3.6]{guolak}.
Let $(X,\| \, \|)$ be a real Banach space, $U\subset X$ be an open bounded set with $0\in U$, 
$K\subset X$ be a cone, $T:K\cap \overline{U}\to K$
be compact and suppose that 
$$\inf_{x\in K\cap \partial U}\|Tx\|>0.$$
Then there exist $\lambda_{0}\in (0,+\infty)$ and $x_{0}\in K\cap \partial U$
such that $x_{0}=\lambda_{0} Tx_{0}$.
\end{thm}
Here, by means of classical fixed point index, we prove a different version of the Birkhoff-Kellogg result, set within the context of \emph{affine} cones. Our result is motivated by the study of retarded functional differential equations. In fact,
when dealing with the solvability of a boundary value problem with delays and initial data, it is somewhat natural to rewrite it in the form of a perturbed integral equation and to seek the solutions of this equation in an affine cone. In particular, the authors in~\cite{acgi16} proved, by means of fixed point index in an affine cone of continuous functions, the existence of multiple nontrivial solutions of the 
perturbed Hammerstein integral equations of the type 
\begin{equation*}\label{eqhamm-ci-old}
u(t)=\psi (t) + \int_{0}^{1} k(t,s)g(s)F(s,u_s)\,ds + \gamma(t) \alpha[u],
\end{equation*}
where $\alpha[\cdot]$ is a \emph{linear} functional in the space $C[0,1]$ given by Stieltjes integral, namely
\begin{equation*}\label{nbc-old}
\alpha[u]=\int_0^1 u(s)\,dA(s).
\end{equation*}
Here we discuss the solvability of the perturbed integral equations
\begin{equation*}\label{eqhamm-intro}
u(t)=\psi(t) + \lambda \Bigl( \int_{0}^{1} k(t,s)g(s)F(s,u_s)\,ds +  \gamma(t) B[u] \Bigr),
\end{equation*}
where $\lambda$ is a non-negative parameter and $B[\cdot]$ is a (not necessarily linear) functional in $C^1([-r, 1], \R)$.
The functional $B[\cdot]$ allows to cover the interesting case of nonlinear and nonlocal boundary conditions (BCs) that can occur in the differential problems; there exists a wide literature on these kind of BCs, we refer the reader to the reviews~\cite{Cabada1, Conti, rma, Picone, sotiris, Stik, Whyburn} and the manuscripts~\cite{Chris-bj, kttmna, jw-gi-jlms}. 
We mention, in particular, the contributions of Mawhin and co-authors in this area of research, see for example~\cite{maw18}.
Note that, in the applications, the functional $B[\cdot]$ can also take into account of the past state of the system.

 As a toy model, we discuss the solvability 
of the 
following class of third order parameter-dependent functional differential equations with functional BCs.
\begin{equation*}
u'''(t)+\lambda F(t,u_t)=0,\ t \in [0,1],
\end{equation*}
with initial conditions
\begin{equation*}\label{kmt-intro}
u(t)=\psi(t),\ t \in [-r,0],
\end{equation*}
and one of the following  BCs
\begin{equation*}\label{fbc1-intro}
u(0)=u'(0)=0,\;  u(1)=\lambda B[u],
\end{equation*}
\begin{equation*}\label{fbc2-intro}
u(0)=u'(0)=0,\;  u'(1)=\lambda B[u],
\end{equation*}
\begin{equation*}\label{fbc3-intro}
u(0)=u'(0)=0,\;  u''(1)= \lambda B[u].
\end{equation*}
Third order functional differential equations with nonlocal boundary terms have been studied in the past, we mention here, for example, the work of Tsamatos~\cite{Tsamatos} and the subsequent papers~\cite{du, yaduge, liyuji}.

As far as we are aware of, our Birkhoff--Kellogg type result (Theorem \ref{BK-translates} below) is new and complements the interesting topological  results in affine cones proved
by Djebali and Mebarki~\cite{djeb2014}.
On the other hand, we also complement the existence results of \cite{acgi16}; this is illustrated in the case of a delay differential equation. In fact here we can deal with equations of the type
 \begin{equation*}\label{eqdelay-intro}
u'''(t)=f(t,u(t),u'(t),u(t-r_1), u'(t-r_2)),\ t \in [0,1]
\end{equation*}
in which we allow the dependence also in the derivative of the solution and we consider the presence of possibly different time-lags.

 \section{Fixed points on translates of a cone}
We require some knowledge of the classical fixed point index for
compact maps, see for example \cite{amannjfa, amann, guolak} for
further information. 
If $\Omega$ is a bounded open subset  (in the relative
topology) of a cone $K$ in a real Banach space we denote by $\overline{\Omega}$ and $\partial \Omega$
the closure and the boundary of $\Omega$ relative to $K$.
Given $y\in X$, we can consider the \emph{translate} of a cone $K$, namely
$$
K_y:=y+K=\{y+x: x\in K\}.
$$
When $D$ is an open
bounded subset of $X$ we write $D_{K_y}=D \cap K_y$, an open subset of
$K_y$.

The following Lemma is a direct consequence of classical results from fixed point index theory (whose properties are analogous to those of the Leray-Schauder degree); a detailed proof can be found, for example, in~\cite{acgi16}. 

\begin{lem} 
Let $(X,\| \, \|)$ be a real Banach space, $K\subset X$ be a cone and
 $D\subset X$ be an open bounded set with $y \in D_{K_y}$ and
$\overline{D}_{K_y}\ne K_y$. Assume that $\F:\overline{D}_{K_y}\to K_y$ is
a compact map such that $x\neq \F x$ for $x\in \partial D_{K_y}$. Then
the fixed point index $i_{K_y}(\F, D_{K_y})$ has the following properties.
\begin{itemize}
\item[(1)] If there
exists $e\in K\setminus \{0\}$ such that $x\neq \F x+\sigma e$ for
all $x\in \partial D_{K_y}$ and all $\sigma >0$, then $i_{K_y}(\F, D_{K_y})=0$.
\item[(2)] If 
$\mu (x-y) \neq \F x-y$
for all $x\in
\partial D_{K_y}$ and for every $\mu \geq 1$, then $i_{K_y}(\F, D_{K_y})=1$.
\item[(3)] Let $D'$ be open in $X$ with
$\overline{D'}\subset D_{K_y}$. If $i_{K_y}(\F, D_{K_y})=1$ and
$i_{K_y}(\F, D'_{K_y})=0$, then $\F$ has a fixed point in
$D_{K_y}\setminus \overline{D'}_{K_y}$. The same result holds if
$i_{K_y}(\F, D_{K_y})=0$ and $i_{K_y}(\F, D'_{K_y})=1$.
\end{itemize}{}
\end{lem}
Our Birkhoff-Kellogg type result is a consequence of the Solution and Homotopy invariance properties of the index. The result reads as follows.
\begin{thm} \label{BK-translates}
Let $(X,\| \, \|)$ be a real Banach space, $K\subset X$ be a cone and
 $D\subset X$ be an open bounded set with $y \in D_{K_y}$ and
$\overline{D}_{K_y}\ne K_y$. Assume that $\F:\overline{D}_{K_y}\to K$ is
a compact map and consider the operator
$$
\F_{(y,\lambda)}:=y+\lambda \F,
$$
where $\lambda \in \R$. Assume that there exists $\bar{\lambda}\in (0,+\infty)$ such that $i_{K_y}(\F_{(y,\bar{\lambda})}, D_{K_y})=0$. Then there exist $x^*\in \partial D_{K_y}$ and $\lambda^*\in (0,\bar{\lambda})$ such that $x^*= y+\lambda^* \F (x^*)$.
\end{thm}
\begin{proof}
First of all note that we have $i_{K_y}(y, D_{K_y})=1$ by the Solution property of the index.
Consider the map
$H:[0,1] \times \overline{D}_{K_y} \to E$ defined by
$H(t,x)= y+t\bar{\lambda} \F (x)$. Note that $H$ is a compact map with values in $K_y$. If there exist $t^*\in (0,1)$ and 
$x\in \partial D_{K_y}$
such that $x= y+t^*\bar{\lambda} \F (x)$ we are done. 
If it does not happen, the fixed point index is defined for $y+t\bar{\lambda} \F$ for every $t\in [0,1]$ and by the Homotopy invariance property we obtain
$$1=i_{K_y}(y, D_{K_y})=i_{K_y}(\F_{(y,\bar{\lambda})}, D_{K_y})=0$$
and the result follows.
\end{proof}
As a Corollary of Theorem~\ref{BK-translates} we exhibit a norm-type Birkhoff-Kellogg-result which can be useful in applications. In order to prove it, we make use of the following proposition.
\begin{pro}[Proposition 2.1 of~\cite{djeb2014}] \label{Proposition2.1}
Let $(X,\| \, \|)$ be a real Banach space, $K\subset X$ be a cone and
 $D\subset X$ be an open bounded set with $y \in D_{K_y}$ and
$\overline{D}_{K_y}\ne K_y$. Assume that $\F:\overline{D}_{K_y}\to K$ is
a compact map and assume that 
\begin{itemize}
\item[(a)] $\displaystyle\inf_{x\in \partial D_{K_y}}\|\F (x)\|>0$
\item[(b)] $\F(x)\neq \mu (x-y)$ for every  $x\in \partial D_{K_y}$ and $\mu \in (0,1]$.
\end{itemize}
Then, $i_{K_y}(\F, D_{K_y})=0$.
\end{pro}
We can now state our norm-type result, which can be seen as an affine version of Theorem~\ref{B-K}.
\begin{cor} \label{BK-transl-norm}
Let $(X,\| \, \|)$ be a real Banach space, $K\subset X$ be a cone and
 $D\subset X$ be an open bounded set with $y \in D_{K_y}$ and
$\overline{D}_{K_y}\ne K_y$. Assume that $\F:\overline{D}_{K_y}\to K$ is
a compact map and assume that 
$$
\inf_{x\in \partial D_{K_y}}\|\F (x)\|>0.
$$
Then there exist $x^*\in \partial D_{K_y}$ and $\lambda^*\in (0,+\infty)$ such that $x^*= y+\lambda^* \F (x^*)$.
\end{cor}
\begin{proof}
We make use of Proposition \ref{Proposition2.1} 
with the map $\bar{\lambda}\F$ in place of $\F$.

We proceed by contradiction and assume that there exist $x_1\in \partial D_{K_y}$ and $\mu_1 \in (0,1]$ such that $\bar{\lambda}\F(x_1)= \mu_1 (x_1-y)$.
Take $R= \sup_{x\in \partial D_{K_y}}\|x\|$, then we have
$$\bar{\lambda} \cdot \inf_{x\in \partial D_{K_y}}\|\F (x)\| \leq \|\bar{\lambda}\F(x_1)\|= \|\mu_1 (x_1-y)\|\leq \| x_1-y\|\leq \| x_1\| + \|y\|\leq R+ \|y\|,$$
a contradiction if
$$
\bar{\lambda}>\frac{R+ \|y\|}{\inf_{x\in \partial D_{K_y}}\|\F (x)\|}.
$$
Then, the result then follows from Theorem~\ref{BK-translates}.
\end{proof}

\section{Positive solutions for a class of perturbed integral equations}
Given a compact interval $I\subset \R$, by $C^1(I, \R)$ we mean the Banach space of the
continuously differentiable functions defined on $I$ with the norm $$\|u\|_{I,1}:=\max\{\|u\|_{I,\infty},\|u'\|_{I,\infty}\},$$ where
$\|u\|_{I,\infty}:=\sup_{t\in I}|u(t)|$.

Given $r>0$ and a continuous function $u: J \to \R$, defined on a real interval
$J$, and given $t \in \R$ such that $[t-r, t] \subseteq J$, we adopt the
standard notation $u_t : [-r, 0] \to \R$ for the function defined by
$u_t (\theta) = u(t + \theta)$.

We consider the following integral equation in the space $C^1([-r, 1], \R)$:
\begin{equation}\label{eqhamm}
u(t)=\psi(t) + \lambda \Bigl( \int_{0}^{1} k(t,s)g(s)F(s,u_s)\,ds +  \gamma(t) B[u] \Bigr)=:\psi(t) + \lambda \F u(t),\quad t \in [-r, 1]
\end{equation}
where $B$ is a suitable (possibly nonlinear) functional in the space $C^1([-r, 1], \R)$.

We require the following assumptions on $r$ as well as on the maps $F$, $k$, $\psi$, $\gamma$  
and $g$ that occur in~\eqref{eqhamm}.
\begin{enumerate}
\item [$(C_{1})$] The function $\psi: [-r,1] \to [0,+\infty)$ is continuously differentiable and
such that $\psi(t)=\psi'(t)=0$ for all $t\in[0,1]$.
\item [$(C_{2})$] The kernel $k:[-r,1] \times [0,1] \to [0,+\infty)$ is measurable,
verifies $k(t,s)=0$ for all $t\in[-r,0]$ and almost every (a.\,e.) $s \in[0,1]$, and
for every $\bar t \in
[0,1]$ we have
\begin{equation*}
\lim_{t \to \bar t} |k(t,s)-k(\bar t,s)|=0 \;\text{ for a.\,e. } s \in
[0,1].
\end{equation*}{}
\item [$(C_{3})$]
 For a.e. $s$, the partial derivative $\partial_{t}k(t,s)$ is continuous in $t$ except at the point $t=s$ where there can be a jump discontinuity, that is, right and left limits both exist, 
and there exists $\Psi \in L^1(0, 1)$ such that $|\partial_{t}k(t,s)| \le \Psi(s)$ for $t \in [0,1]$ and a.e. $s \in [0,1]$.
{}
\item [ $(C_{4})$] The function $g:[0,1] \to \R$ is measurable, $g(t) \geq 0$ a.\,e. $t \in [0,1]$,
and satisfies that $g\,\Phi \in L^1[0,1]$
and $\int_a^b \Phi(s)g(s)\,ds >0$.{}
\item [$(C_{5})$] $F: [0,1] \times C^1([-r, 0], \R) \to [0,\infty)$ is an operator that satisfies some 
Carath\'eodory-type conditions (see also \cite{halelunel}); namely,
for each $\phi$, $t \mapsto F(t,\phi)$ is measurable and for a.\,e. $t$, $\phi \mapsto F(t,\phi)$ is continuous. Furthermore,
for each $R>0$, there exists $\varphi_{R} \in
L^{\infty}[0,1]$ such that{}
\begin{equation*}
F(t,\phi) \le \varphi_{R}(t) \ \text{for all} \ \phi \in C^1([-r, 0], \R)
\ \text{with} \ \|\phi\|_{[-r,0],1} \le R,\ \text{and a.\,e.}\ t\in [0,1].
\end{equation*}{}
\item [ $(C_{6})$]
The function $\gamma: [-r,1] \to [0,\infty)$ is continuous differentiable, and such that $\gamma(t)=\gamma'(t)=0$ for all $t\in[-r,0]$.
\end{enumerate}

In the Banach space $C^1([-r, 1], \R)$ we define the cone of non-negative functions
$$
K_0=\{u\in C^1([-r, 1], \R): u(t)\geq 0\ \text{for every}\ t\in[-r,1]\ \text{and}\ u(t)=u'(t)= 0 \ \text{for every}\ t\in[-r,0]\}.
$$
Note that the function 
$$w(t)=\begin{cases}
0,\ & t \in [-r,0], \\
t^2,\ & t \in [0,1],
\end{cases}
$$
belongs to $K_0$, hence $K_0 \neq \{0\}$.

We consider the following translate of the cone $K_0$,
$$K_\psi=\psi + K_0 = \{\psi +u : u \in K_0\}.$$
\begin{defn}
We define the following subsets of $C^1([-r, 1], \R)$: $$K_{0,\rho}:=\{u\in K_0: \|u\|_{[0,1],1} <\rho\},\quad 
K_{\psi,\rho}:= \psi + K_{0,\rho}.$$
\end{defn}

The following theorem provides an existence result for equation~\eqref{eqhamm}: here we obtain a non-trivial solution within the cone $K_\psi$ with fixed norm and a corresponding positive parameter.
\begin{thm}\label{eigen}Let $\rho \in (0,+\infty)$ and assume the following further conditions hold. 
\begin{itemize}

\item[$(a)$] 
There exist $\underline{\delta}_{\rho} \in C([0,1],\R_+)$ such that
\begin{equation*}
F(t,\phi)\ge \underline{\delta}_{\rho}(t),\ \text{for every}\ (t,\phi)\in [0,1] \times \partial K_{\psi,\rho}.
\end{equation*}

\item[$(b)$] 
$B: \overline K_{\psi,\rho} \to \R_{+}$ is continuous and bounded. Let $\underline{\eta}_{\rho}\in [0,+\infty)$ be such that 
\begin{equation*}
B[u]\geq \underline{\eta}_{\rho},\ \text{for every}\ u\in \partial K_{\psi,\rho}.
\end{equation*}
\item[(c)] 
The inequality 
\begin{equation}\label{condc}
\sup_{t\in [0,1]}\Bigl\{ \gamma(t)\underline{\eta}_{\rho}+\int_{0}^{1}  k(t,s)g(s) \underline{\delta}_{\rho} (s)\,ds\Bigr\}>0 
\end{equation}
holds.
\end{itemize}
Then there exist $\lambda_\rho$ and $u_{\rho}\in \partial K_{\psi,\rho}$ such that the 
 integral equation~\eqref{eqhamm} 
 is satisfied.
\end{thm}

\begin{proof}
Consider the operator $\F u$ defined in~\eqref{eqhamm}. 
Due to the assumptions above, $\F$ maps
$\overline{K}_{\psi,\rho}$ into $K_0$ and is compact. The compactness of  the Hammerstein integral operator is a consequence of the regularity assumptions on the terms occurring in it combined with a careful use of  the Arzel\`{a}-Ascoli theorem (see~\cite{Webb-Cpt}), while the perturbation 
$\gamma(t) B[ \cdot ]$ is a finite rank operator. 

Take $u\in \partial K_{\psi, \rho}$, then we have

\begin{multline}\label{lwest}
\|  \F u\|_{[-r, 1],1}\geq   \|  \F u\|_{[-r, 1],\infty}=\sup_{t\in [0, 1]} \Bigl| \int_{0}^{1} k(t,s)g(s)F(s,u_s)\,ds +  \gamma(t) B[u]\Bigr |\\
\geq \sup_{t\in [0,1]}\Bigl\{ \gamma(t)\underline{\eta}_{\rho}+\int_{0}^{1}  k(t,s)g(s) \underline{\delta}_{\rho} (s)\,ds\Bigr\}.
\end{multline}

Note that the RHS of \eqref{lwest} does not depend on the particular $u$ chosen. Therefore we have
$$
\inf_{u\in \partial K_{\psi, \rho}}\|  \F u\|_{[-r, 1],1} \geq \sup_{t\in [0,1]}\Bigl\{ \gamma(t)\underline{\eta}_{\rho}+\int_{0}^{1}  k(t,s)g(s) \underline{\delta}_{\rho} (s)\,ds\Bigr\}>0,
$$
and the result follows by Corollary~\ref{BK-transl-norm}.
\end{proof}

\section{an application}
We now apply the previous results to the following class of third order functional differential equations with functional BCs.
\begin{equation}\label{FDE}
u'''(t)+\lambda F(t,u_t)=0,\ t \in [0,1],
\end{equation}
with initial conditions
\begin{equation}\label{kmt-i}
u(t)=\psi(t),\ t \in [-r,0],
\end{equation}
and one of the following boundary conditions (BCs)
\begin{equation}\label{fbc1}
u(0)=u'(0)=0,\;  u(1)=\lambda B[u],
\end{equation}
\begin{equation}\label{fbc2}
u(0)=u'(0)=0,\;  u'(1)=\lambda B[u],
\end{equation}
\begin{equation}\label{fbc3}
u(0)=u'(0)=0,\;  u''(1)= \lambda B[u].
\end{equation}
We begin by considering some auxiliary problems. 

First of all note that the solution of the ODE $-u'''=y$ under the BCs
\begin{equation}\label{fbc1h}
u(0)=u'(0)=u(1)=0,
\end{equation}
\begin{equation}\label{fbc2h}
u(0)=u'(0)=u'(1)=0,
\end{equation}
\begin{equation}\label{fbc3h}
u(0)=u'(0)=u''(1)=0,
\end{equation}
in the interval $[0,1]$ is given by
$$
u(t)= \int_{0}^{1}\hat{k}_i(t,s)y(s)ds,
$$
where the Green's function is
$$
\hat{k}_1(t,s)= \frac{1}{2}\begin{cases}s(1-t)(2t-ts-s),\ &s\leq t,\\ (1-s)^2t^2,\ &s\geq t,
\end{cases}
$$
in the case of the BCs~\eqref{fbc1},
$$
\hat{k}_2(t,s)= \frac{1}{2}\begin{cases}(2t-t^2-s)s,\ &s\leq t,\\ (1-s)t^2,\ &s\geq t,
\end{cases}
$$
for the BCs~\eqref{fbc2} and 
$$
\hat{k}_3(t,s)= \frac{1}{2}\begin{cases}s(2t-s),\ &s\leq t,\\ t^2,\ &s\geq t,
\end{cases}
$$
for the BCs~\eqref{fbc3}. Furthermore note that the function
$$
\hat{\gamma}_1(t):=t^2
$$ is the unique solution of the BVP
\begin{equation*}
\hat{\gamma}'''(t)=0,\ \hat{\gamma}(0)=\hat{\gamma}'(0)=0,\; \hat{\gamma}(1)=1,
\end{equation*}
while the functions
$$
\hat{\gamma}_2(t)\equiv\hat{\gamma}_3(t):=\frac{1}{2}t^2
$$ solve the BVPs
\begin{equation*}
\hat{\gamma}'''(t)=0,\ \hat{\gamma}(0)=\hat{\gamma}'(0)=0,\; \hat{\gamma}'(1)=1.
\end{equation*}
\begin{equation*}
\hat{\gamma}'''(t)=0,\ \hat{\gamma}(0)=\hat{\gamma}'(0)=0,\; \hat{\gamma}''(1)=1.
\end{equation*}
By routine calculations (see also \cite{guen-had-ben-2018, sun-li-2010}) one obtains the following proposition.
\begin{pro}
For every $i=1,2,3$, we have:
\begin{enumerate}
\item $\hat{k}_i$ is continuous and non-negative in $[0,1]\times [0,1]$ and the partial derivative $\partial_{t}k(t,s)$ is continuous in $t\in [0,1]$ for every $s \in [0,1]$.
\item $\hat{\gamma}_i$ is non-negative and continuously differentiable in $[0,1]$.
\end{enumerate}
\end{pro}
Due to the above setting, the functional boundary value problem (FBVP) \eqref{FDE}-\eqref{kmt-i}-\eqref{fbc1} can  be rewritten in the form \eqref{eqhamm},
where $\gamma_1(t):=H(t) \hat{\gamma}_1(t)$ and $k_1(t,s):=H(t) \hat{k}_1(t,s)$ with
$$
H(\tau)=
\begin{cases}
1,\ & \tau \geq 0, \\0,\ & \tau<0,
\end{cases}
$$
and, provided that $\psi, F, B$ possess a suitable behaviour, Theorem \ref{eigen} can be applied directly; this fact holds also in the case of the FBVPs \eqref{FDE}-\eqref{kmt-i}-\eqref{fbc2} and \eqref{FDE}-\eqref{kmt-i}-\eqref{fbc3}.

We now describe the applicability of our theory to the context of \emph{delay differential equations}. Namely, 
 let $f:[0,1]\times\R_+\times\R \times\R_+\times\R \to [0,\infty)$ be a given Carath\'eodory map, and consider the equation
 \begin{equation}\label{eqdelay}
u'''(t)=f(t,u(t),u'(t),u(t-r_1), u'(t-r_2)),\ t \in [0,1],
\end{equation}
where $r_1$ and $r_2$ are positive and fixed (possibly different).
We can apply the techniques developed in this paper to the equation 
\eqref{eqdelay} with initial condition \eqref{kmt-i} along with one of the BCs \eqref{fbc1}, \eqref{fbc2}, \eqref{fbc3}.
To see this, observe that \eqref{eqdelay} is a special case of the functional equation \eqref{FDE}, in which taking $r:= \max\{r_1,r_2\}$,
 the operator $F: [0,1] \times C^1([-r, 0], \R) \to [0,\infty)$ is defined by
$$
F(t,\phi)=f(t,\phi(0),\phi'(0),\phi(-r_1),\phi'(-r_2)).
$$ 
Such an operator satisfies the above condition $(C_{5})$ provided that  the following assumption on the map $f$ is verified:

[$(C'_{5})$] For each $R>0$, there exists $\varphi^*_{R} \in
L^{\infty}[0,1]$ such that{}
\begin{align*}
f(t,u,v,p,q) \le \varphi^*_{R}(t) \;\text{ for all } \; (u,v,p,q) \in \R_+\times\R \times\R_+\times\R \\
\;\text{ with } \; 0\le u,p \le R,\;|v| \le R,\;|q| \le R,\;\text{ and a.\,e. } \; t\in [0,1].
\end{align*}

To better illustrate the growth conditions we now provide a specific example.

\begin{ex}
We adapt the nonlinearities studied in Example 2.6 of~\cite{giaml22} to the context of delay equations by consider the family of FBVPs
\begin{equation}\label{eqdiffex}
u'''(t)+\lambda te^{u(t)+(u'(t-\frac12))^2}(1+(u'(t))^2+(u(t-\frac13))^2),\ t \in (0,1),
\end{equation}
with the initial condition
\begin{equation}\label{inex}
u(t)=\psi(t), t \in [-\frac12,0],
\end{equation}
with $\psi(t)=H(-t)t^2$,
and one of the three BCs~\eqref{fbc1h}, \eqref{fbc2h}, \eqref{fbc3h},
where we fix
$$
B[u]=\lambda \Bigl(\frac{1}{1+(u(\frac{1}{2}))^2}+\int_{-\frac12}^1t^3(u'(t))^2\, dt \Bigr).
$$
Now choose $\rho\in (0,+\infty)$. Thus we may take 
$$\underline{\eta}_{\rho}(t)=\frac{1}{1+\rho^2}, \underline{\delta}_{\rho}(t)=t.$$ 

Therefore, for every $i=1,2,3$, we have
$$
\sup_{t\in [0,1]}\Bigl\{ \frac{\gamma_i(t) }{1+\rho^2}+\int_{0}^{1}  k_i(t,s)t\,ds\Bigr\}\geq \frac{1}{2(1+\rho^2)} > 0, 
$$
which implies that \eqref{condc} is satisfied for every $\rho \in (0,+\infty)$.

Thus we can apply Theorem~\ref{eigen}, obtaining uncountably many pairs of solutions and parameters $(u_{\rho}, \lambda_{\rho})$ for the FBVPs~\eqref{eqdiffex}-\eqref{inex}-\eqref{fbc1h}, \eqref{eqdiffex}-\eqref{inex}-\eqref{fbc2h} and \eqref{eqdiffex}-\eqref{inex}-\eqref{fbc3h}.
\end{ex}

\section*{Acknowledgements}
The authors were partially supported by
 the Gruppo Nazionale per l'Analisi Matematica, la Probabilit\`a e le loro Applicazioni (GNAMPA) of the Istituto Nazionale di Alta Matematica (INdAM).
G.~Infante is a member of the UMI Group TAA “Approximation Theory and Applications”.

\end{document}